\documentclass[12pt,a4]{article}
\usepackage{amsmath,amsfonts,amsthm,amstext}
\usepackage{epic}   %from Kuo
\usepackage{eepic}  %from Kuo
\usepackage{amssymb,latexsym}
\usepackage{epsfig} %from cocke
\usepackage{graphicx} %from nuno

\newskip\foobarskip
\newlength{\pagewidth}
\newtheorem{prop}{Proposition}[section]
\newtheorem{theo}[prop]{Theorem}

\newtheorem{ex}[prop]{Example}
\newtheorem{defn}[prop]{Definition}
\newtheorem{rem}[prop]{Remark}
\newtheorem{lem}[prop]{Lemma}

\newcommand{\R}{\mathbb{R}}
\newcommand{\C}{\mathbb{C}}

\newcommand{\Z}{\mathbb{Z}}

\newenvironment{prf}
 {\begin{trivlist} \item[\hskip \labelsep {\bf Proof}\hspace*{3 mm}]}
 {\hfill$\Box$\end{trivlist}}

\begin{document}
  \title{Counting isolated singularities in germs of applications
 ${\C^n,0\to \C^p,0}$\,\,\,\,$n<p$}
  \author{V. H. Jorge  P\'erez\thanks{Work partially supported CAPES-PROCAD}}
  \date{}
  \maketitle

\begin{abstract}

In this paper we  give a formula for counting the number of
isolated stable singularities of a stable perturbation of corank 1
germs $f:\C^n,0\to \C^p,0$ with $n<p$ that appear in the image
$f(\C^n).$
 We also define a set of
${\cal A}$-invariants  and show that their finiteness  is a
necessary and sufficient condition for the  ${\cal A}$-finiteness
of the germ $f.$

\end{abstract}

  \section{Introduction}
  Let $f:{\C^n,0\to \C^p,0},$ be a germ of  a map with a finite number of
  isolated stable singularities or zero-dimensional stable
types  in the discriminant of $f$. For example, if $n=2$, $p=3$,
the isolated stable singularities in the hypersurface $f(\C^2)$
are the cross-caps  and triple points. {\it The main problem is to
find an algebraic formula  for counting the number of isolated
singularities of  $f$.}

When $n=2$, $p=3$, D. Mond  shows in \cite{mond} that the number
of cross-cap ($C(f)$) and the number of triple points ($T(f)$) are
given by the dimensions of local algebras associated to $f$.

 When $n=2$, $p=2$ J. Rieger \cite{rie} shows that the
number of cusps is given by the dimension of a local algebra
associate to $f$ in the case where $f$ is of corank 1. T. Gaffney
and D. Mond  \cite{GM} give formulae   for both the number of
cusps and the number of double points for a general
finitely-determined map-germ $\C^2,0\to \C^2,0$.

 When $n=p$ W. Marar, J. Montaldi and M.Ruas \cite{ruas}
give formulae for calculating the  isolated stable singularities
associated to $f$ in the case when $f$ is weighted homogeneous and
of corank 1.

In this work we consider the analogous problem for map-germs
$f:\C^n,0\to \C^p,0$ with $n<p$  and give a formula for counting
the number of all isolated stable singularities. In particular, if
$f$ is weighted homogeneous,  we give a  formula for these numbers
in terms of the weights  of the variables and the degrees of each
component of $f$.

We also define a set of ${\cal A}$-invariants  and show that their
finiteness  is a necessary and sufficient conditions for the
${\cal A}$-finiteness of the germ $f$.

\section{Stable types}
Our notation are standard  in singularity theory. We denote by
${\cal A}$ the group $Diff(\C^n,0)\times Diff(\C^p,0);$ this acts
on ${\cal O}(n,p)$ the space of germs $\C^n,0 \to \C^p,0$, by
composition on the right and on the left.

A $d$-parameter {\it unfolding} of a map-germ $f_0\in {\cal
O}(n,p)$ is a germ $F\in {\cal O}(n+d,p+d)$ of the form
$F(x,u)=(f(x,u),u),$ with $f(x,0)=f_0$. A $c$-parameter unfolding
$F'$ de $f_0$ is {\it induced} from a $d$-parameter unfolding $F$
by a germ $h:\C^c,0 \to \C^d,0$ if $(f'(x,v),v)=(f(x,h(v)),v)$, An
unfolding $F$ of $f_0$ is ${\cal A}$-{\it versal} if every other
unfolding of $f_0$ is ${\cal A}$-equivalent to an unfolding
induced from $F$. An ${\cal A}$-versal unfolding of $f_0$
contains, up to ${\cal A}$-equivalence, every other unfolding of
$f_0.$

A map-germ $f:\C^n,0 \to \C^p,0$ is $k$-${\cal A}$-{\it
determined} if, whenever $g:\C^n,0 \to \C^p,0$ and
$j^kf(0)=j^kg(0)$, $g$ is ${\cal A}$-equivalent to $f$, and ${\cal
A}$-{\it finite} if it is $k$-${\cal A}$- determined for some
$k<\infty$. In this paper we shall refer almost exclusively to
${\cal A}$-finite. We say $f$ is a {\it stable germ}  if every
nearby germ is ${\cal A}$-equivalent to $f$.

Let $F$ be versal unfolding of $f$. We say that a
 stable type ${\cal Q}$ appears in $F$ if, for any
representative $F=(id,f_s)$ of $F$, there exists a point $(s,y)\in
{\C}^s\times {\C}^p$ such that the germ $f_s:{\C}^n,S\to {\C}^p,y$
is a stable germ of type $\cal Q $, where $S=f^{-1}(y)\cap
\Sigma(f_s).$ We call $(s,y)$ and the points $(s,x)$ with $x\in S$
points of stable type $\cal Q$ in the target and in the source,
respectively. If $f$ is stable, we denote the set of points in
$\C^s\times\C^p$ of type $\cal Q$ by ${\cal Q}(f)$ and set ${\cal
Q}_S(f)=f^{-1}({\cal Q}(f))-{\cal Q}_{\Sigma}(f),$ where ${\cal
Q}_{\Sigma}(f)$ denotes $f^{-1}({\cal Q}(f)) \cap \Sigma(f).$

If $f$ is ${\cal A}$-finite, we denote $\overline{{\cal Q}(f)} =
(\{0 \} \times {\C}^p)\cap {\overline{{\cal Q}(F)}}$ and
$\overline{{\cal Q}_S(f)}=(\{0 \} \times {\C}^n)\cap
{\overline{{\cal Q}_S(F)}},$ $\overline{{\cal Q}_{\Sigma}(f)}=(\{0
\} \times {\C}^n)\cap {\overline{{\cal Q}_{\Sigma}(F)}}.$

We say that $\cal Q$ is a {\it zero-dimensional stable type or
isolated stable singularity} for the pair $(n,p)$ if ${\cal Q}(f)$
has dimension $0$ where $f$ is a representative of the stable type
$\cal Q$. We observe that ${\cal Q}(F)$ is a close analytic set
since $\overline{{\cal Q}(F)}=\cap F(j^{(p+1)
}F^{-1}(\overline{{\cal A}{z_i}})),$ where $z_i$ is the $p+1$ jet
of the stable type ${\cal Q}$ and ${\cal A}{z_i} $ is the ${\cal
A}-$orbit of $z_i $.

An ${\cal A}$-finite  germ  $f$ has a discrete stable type if
there exist a versal unfolding $F$ of $f$ in which  only a finite
number of stable types appear. An ${\cal A}$-finite germ  $f$ has
a discrete stable type if the pair $(n,p)$ is in the nice
dimensions (\cite{mater}).

In the next sections we define explicitly the zero-dimensional
stable types
 or stable isolated singularities, that is, we
 define  the ideal that defines each zero dimensional stable
type  of $f$.

%\rm{
%\begin{ex}
%For ${\cal A}-$finite map germs $f:\C^2\to \C^3$ of corank 1, the
%types stable singularities  of $f(\C^2)$ are:
%$D_1^2(f,(1,1))=D^2(f)$ is the set of  points  double,
%$f(D_1^1(f,(2)))$ is the set of points cross-cap and
%$D^3(f,(1,1,1))=D^3(f)$ is the set of points triples.
%\end{ex}
%}
%\medskip
%
%For ${\cal A}-$finite map germs $f\in {\cal O}(3,3)$ of corank 1,
%os tipos estaveis s\~ao $D_1^1(f,(1))=\Sigma(f)$  is the singular
%set of $f$, $D_1^2(f,(1,1))=D^2(f)$ is the set of double points
%of $f$, $D_1^1(f,(2))=\Sigma^{1,1}(f)$ is a cuspidal curve,
%$D^3(f,(1,1,1))=D^3(f)$ is the set of triple points and
%$D_1^1(f,(3))= \Sigma^{1,1,1}(f)$ is the Swallowtail of $f$.
%\end{ex}

\section{Stable isolated singularities}
Consider a versal unfolding  $F$ of $f$ with base $\C^d,0$,
 $$
 \begin{array}{c}
 F:\C^{n-1}\times\C\times\C^d,0\to \C^p\times\C^d,0\\
F(x,y,u)=(x,h_1(x,y,u),...,h_{p-n+1}(x,y,u),u)
\end{array}
$$

For each $f$ we associate a partition ${\cal P}
=(r_1,...,r_{\ell})$ of $k$ where $p-k(p-n+1)+\ell=0$ and consider
$\tilde{V}({\cal P})\subset \C^{n-1}\times\C^{\ell}\times\C^d,0$
defined by

$$
\tilde{V}({\cal P})=clos\left\{\begin{array}{l}(x,{\bf z},u)\in
\C^{n-1}\times \C^{\ell}\times\C^d |\end{array} \begin{array}{l}
{\bf z}=(z_1,\cdots,z_{\ell}),\,\,\,\,
z_i\neq z_j\\
 F(x,z_i,u)=F(x,z_j,u)\\
F_u\,\, \mbox{
has  a}\,\,  {\cal Q}\,\,\mbox{0-dimensinal }\\
\mbox{stable singularity at} (x,z_j)
\end{array}\right\}
$$

\noindent where `clos' means the analytic closure in
$\C^{n-1}\times \C^{\ell}\times\C^d$. The varieties
 $\tilde{V}({\cal P})$ are called zero-schemes and are
related to the $0$-stable types (see \cite{gaff}).

Let $\pi_{\cal P}:\tilde{V}({\cal P})\to \C^d$ be the restriction
to $\tilde{V}({\cal P})$ of the cartesian projection
$\C^{n-1}\times \C^{\ell}\times\C^d\to\C^d$. For a generic
$u\in\C^d$, the fiber $\pi^{-1}_{\cal P} (u)$ consists of the
multiple points  where $F_u$ has a ${\cal Q} (f,{\cal P})$
multi-germ that defines a zero dimensional stable type.  We are
thus interested in the degree of $\pi_{\cal P} $.

Let ${{\cal P}}=(r_1,\cdots,r_{\ell})$ be a partition of $m\leq n$
with $r_1\geq r_2\geq \cdots r_{\ell}\geq 1$. Define $N({\cal P})$
to be the order of the sub group of $S_{\ell}$ which fixes ${\cal
P}$. Here $S_{\ell}$ acts on $\R^{\ell}$ by permuting the
coordinates. For example, if ${\cal P}=(4,4,4,2,2,2,1,1)$ we have
$N({\cal P})=(3!)^22!$.

\begin{prop}\label{2} If ${\cal P} =(r_1,...,r_{\ell})$ is
a partition of $k$, then
$$\sharp {\cal Q}(f,{\cal P})=\frac{1}{N({\cal P})}{\deg}(\pi_{\cal P}).$$
\end{prop}
\begin{prf} Let ${\bf y} =({\bf x} ,y_1,...,y_{\ell},u)\in
\tilde{V}({\cal P}) $ and $\sigma \in S_{\ell}$. We have\\
${\bf y_{\sigma}}= (y_{\sigma(1)},...,y_{\sigma({\ell})},u)\in
\tilde{V}({\cal P})$ if and only if $r_{\sigma(j)}=r_j$ for each
$j=1,...,\ell.$ There are $N({\cal P})$ such $\sigma$. The points
${\bf y}$ and ${\bf y_{\sigma}}$ are distinct, but the
corresponding sets $\{{\bf x} ,y_1,...,y_{\ell}\}$  are the same,
and it is the contribution of the multiple points that are conted
in $\sharp{\cal Q }(f,{\cal P} )$.
\end{prf}

Let $\tilde{I}(\tilde{f},{\cal P})$ be the ideal in ${\cal
O}_{n-1+\ell+d}$ that defines  $\tilde{V}({\cal P} )$, and let
$$
I(f,{\cal P})=\frac{(\tilde{I}(\tilde{f},{\cal P})+
<u_1,...,u_d>)}{<u_1,...,u_d>}\subset{\cal O}_{n-1+\ell}
$$
be the ideal corresponding to the intersection of $\tilde{V}({\cal
P})$ with $\C^{n-1+\ell}\times \{0\}.$

It follows from the definition of ${\tilde{I}({{\tilde f},\cal P}
})$ that, at generic points of $\tilde{V}({\cal P})$ with $z_i\neq
z_j$, $ \tilde{I}(\tilde{f},{\cal P})= (\{\frac{\partial^s
  \tilde{f}_j}{\partial z^s}\circ \pi_i({{\cal P}})| j=1,...,p-n+1,1\leq
  s\leq r_i-1,1\leq i \leq m \})+\\
  \,\,\,\,\,\,\,\,\,\,\,\,\,(\{\tilde{f}_j\circ \pi_1({{\cal P}})-\tilde{f}_j\circ
 \pi_i({{\cal P}})|j=1,....p-n+1, 2\leq i\leq m\})$ in ${{\cal
  O}_{{n-1+\ell+d}},(x,z)}$, where $\pi_i({{\cal P}}):\C^{n-1+\ell}\to \C^n,$   $1\leq
  i\leq m,$ are given by $\pi_i({{\cal
  P}})(x,z_1,...,z_\ell)=(x,z_{r_1+...+r_{i-1}+1}).$
In particular, at generic point of $\tilde{V}({\cal P})$, we have
$$ \frac{{\cal
O} _{n-1+\ell+d}}{(m_d,\tilde{I}(\tilde{f},{\cal P}))}\cong
\frac{{\cal O} _{n-1+\ell}}{I(f,{{\cal P} })},$$ where $m_d$ is
the maximal ideal of $\C^d $.

\begin{prop}\label{3} Suppose that $\tilde{V}({\cal P})$
is non-empty. Then
\begin{enumerate}
\item $\tilde{V}({\cal P})$
is smooth of dimension $d$;
\item  $\pi_{\cal P}:\tilde{V}({\cal P} )
\to \C^d$ is finite and $\pi^{-1}_{\cal P} (\pi_{\cal P}
(0))=\{0\};$
 \item The degree of $\pi_{\cal
P}$ coincides with $\dim_{\C}\frac{{\cal O} _{n-1+\ell}}{I(f,{\cal
P})}.$
\end{enumerate}
\end{prop}
\begin{prf}{ 1.} Since $F$ is a versal unfolding,  it is stable,
and the proof follows by Propositions 2.13 in  \cite{marar}.

{ 2.} The projection $\pi_{\cal P}:\tilde{V}({\cal P}) \to \C^d$
is finite. In fact, for $u\in\C^d$ generic, the fiber
$\pi^{-1}_{\cal P} (u)$ is finite and consists of those
multi-points where $f_u$ has a ${\cal Q}(f_u,{\cal P})$ multi
germ. The germ $f_0=f$ is ${\cal A}$-finite, so
 using the geometric criterion  of
Mather-Gaffney (\cite{cg},  \cite{ctc}), it is stable away from
zero. Thus $\pi^{-1}_{\cal P} (\pi_{\cal P} (0))=\{0\}$.

 { 3.} Since $\tilde{V}({\cal P})$ is smooth
 and it is Cohen-Macaulay at zero, the degree of
 $\pi_{\cal P}$ coincides with $\dim_{\C}\frac{{\cal
O} _{n-1+\ell}}{I(f,{\cal P})}$ (see Proposition $5.12$ in
\cite{looig}).
\end{prf}

Propositions \ref{2} and \ref{3} (3) give a formula for computing
 the multiplicities of ${\cal Q}(f_u,{\cal P})$ even in the case
 when
 $f$ is not weighted homogeneous. We have
 $$\sharp {\cal Q}(f,{\cal P})=
 \frac{1}{N({\cal P})}\dim_{\C}\frac{{\cal O} _{n-1+\ell}}{I(f,{\cal
P})}.$$

The $\dim_{\C}\frac{{\cal O} _{n-1+\ell}}{I(f,{\cal P})}$  not is
difficult to calculate, when $I(f,{\cal P})$ can be computed. The
calculation can be done using computer algebra package such as
Singular \cite{singular} or Macaulay.

\begin{ex} \label{exe3} {\rm Let $f:\C^2,0\to\C^3,0$ be the  ${\cal A}$-
finite germ  given by
$$f(x,y)=(x,x^2y+y^6+y^7,xy^2+y^4+y^6+y^9).$$ We chooses the
partition ${\cal P}$ of $k$ such that $k=\frac{3+\ell}{2}\in
\Z^+$. Then ${\cal P}$ has two elements  $(1,1,1)$ and $(2)$. We
have $I(f,(2))=(x^2+6y^5+7y^6,2xy+4y^3+6y^5+9y^8)$ and
$I(f,(1,1,1)) =( -5z_1^6+x^2-4z_1^2z_3^3-4z_1^2z_3^4-2z_1z_3^4-
8z_1^4z_3-6z_1^3z_3^3-8z_1^4z_3^2-4z_1^5-6z_1^3z_3^2-2z_1z_3^5-
10z_3z_1^5,
 3z_1^2z_3^2+3z_1^2z_3^3+2z_1z_3^4+5z_1^4z_3+2z_1z_3^3+
 4z_3z_1^3+6z_1^5+4z_1^3z_3^2+z_3^5+5z_1^4+z_3^4,
-2z_1^3-7z_1^8-6z_3^5z_1^3-10z_3^3z_1^5-8z_3^4z_1^4-
4z_1^2z_3^3-2z_1z_3^4-8z_1^4z_3-4z_1^2z_3-2z_1z_3^2-
2z_1z_3^7-14z_3z_1^7-4z_1^5-6z_1^3z_3^2-4z_1^2z_3^6-12z_3^2z_1^6,
x+8z_1^7+z_3^7+3z_1^2z_3^2+2z_1z_3^6+7z_3z_1^6+2z_1z_3^3+4z_3z_1^3+
6z_1^5z_3^2+2z_1z_3+
3z_3^5z_1^2+4z_3^4z_1^3+5z_3^3z_1^4+z_3^2+3z_1^2+5z_1^4+z_3^4).$

%$i=-5*z^6+x^2-4*z^2*w^3-4*z^2*w^4-2*z*w^4-8*z^4*w-6*z^3*w^3-8*z^4*w^2-
%4*z^5-6*z^3*w^2-2*z*w^5-10*w*z^5,
%3*z^2*w^2+3*z^2*w^3+2*z*w^4+5*z^4*w+2*z*w^3+4*w*z^3+6*z^5+4*z^3*w^2+
%w^5+5*z^4+w^4,
%-2*z^3-7*z^8-6*w^5*z^3-10*w^3*z^5-8*w^4*z^4-4*z^2*w^3-2*z*w^4-8*z^4*w-4*z^2*w-2*z*w^2-
%2*z*w^7-14*w*z^7-4*z^5-6*z^3*w^2-4*z^2*w^6-12*w^2*z^6,
%x+8*z^7+w^7+3*z^2*w^2+2*z*w^6+7*w*z^6+
%2*z*w^3+4*w*z^3+6*z^5*w^2+2*z*w+3*w^5*z^2+4*w^4*z^3+5*w^3*z^4+w^2+3*z^2+5*z^4+w^4$
Using Theorem \ref{3} and  Singular we have $\sharp{\cal Q
}(f,{(2)})=6$  and $\sharp{\cal Q}(f,(1,1,1))=14$.}
\end{ex}

\section{Multiple points}\label{sec4}
We need the following
  definition. Given a continuous mapping $f:X\to Y$ between
  analytic spaces, we define the $k^{tk}$ multiple point space of
  $f$ as $$ D^k(f)=\mbox{clos}\{(x_1,x_2,...,x_k)\in X^k:
  f(x_1)=...=f(x_k)\,\, \mbox{for}\,\, x_i\neq x_j, i\neq j\}. $$
  Suppose $f:\C^n,0\to \C^p,0$ is of corank 1 and is given in the
  form
  $$f(x_1,...,x_{n-1},z)=(x_1,...,x_{n-1},h_1(x,z),...,h_{p-n+1}(x,z)).$$
  If $g:\C^n,0\to \C,0$ is an analytic function then we define $V_i^k(g):\C^{n+k-1},0\to
  \C,0$ to be

  \[
  \left|{\begin{array}{cccccccc}
  1&z_1&\cdots&z_1^{i-1}&g(x,z_1)&z_1^{i+1}&\cdots&z_1^{k-1}\\
  \vdots&\vdots&&\vdots&\vdots&\vdots&&\vdots\\
  1&z_k&\cdots&z_k^{i-1}&g(x,z_k)&z_k^{i+1}&\cdots&z_k^{k-1}
  \end{array}
  }\right|
  {\Huge{/}}
  \left|{
  \begin{array}{cccc}
  1&z_1&\cdots&z_1^{k-1}\\
  \vdots&\vdots&&\vdots\\
  1&z_k&\cdots&z_k^{k-1}.
  \end{array}}\right|
  \]

  \begin{theo}(\cite{marar})\label{teo}
  $D^{k}(f,{\cal P})=D^{k}(f)$ is defined in $\C^{n+k-1}$
  by the ideal ${{\cal I}}^{k}(f)$
  generated by $V_i^{\ell}(h_j(x,z))$
  for all $i=1,...,k-1$, ${\cal P}=(1,...,1)$-$k$
  times and
  $j=1,...,p-n+1$.
   \end{theo}

  In what  follows we take coordinates in
  $\C^{n+k-1}=\C^{n-1}\times \C^k$ to be
  $(x,z)=(x_1,...,x_{n-1},z_1,...,z_k).$

\begin{ex} {\rm For a corank $1$ map-germ $f:\C^n,0\to \C^p,0,$
$D^2(f)$ is defined by $\frac{h_i(x,z_1)-h_i(x,z_2)}{z_1-z_2}$,
$i=1,...,p-n+1.$}
\end{ex}

\begin{defn} Let ${{\cal P}} =(r_1, r_2,...,r_{\ell})$ be a partition of
$k$ where $p-k(p-n+1)+\ell=0$,
 that is, $r_1+r_2+...+r_{\ell}=k$. Let ${{\cal I}}({{\cal P}})$ be the ideal
in ${{\cal O}_{n-1+k}}$ generated by the $k-\ell$ elements
$z_i-z_{i+1}$ for $r_1+r_2+...+r_{j-1}+1\leq i\leq
r_1+r_2+...+r_j-1$, $1\leq j\leq \ell$, and let $\Delta({{\cal
P}})=V({{\cal I}}({{\cal P}}))$. Define $${\cal
J}_{\Delta}(f,{{\cal P}})={\cal  I}^k(f)+{{\cal I}}({{\cal P}})\,
{\mbox{and}}\,\,\, D^{\ell}(f,{{\cal P}})=V({\cal
J}_{\Delta}(f,{{\cal P}})),$$ equipped with the sheaf structure in
${{\cal O}_{n-1+k}}/{\cal J}_{\Delta}(f,{{\cal P}})$.
\end{defn}
\begin{ex}
{\rm For an ${\cal A}-$finite of corank 1 map germ $f\in {\cal
O}(2,3)$, the stable types  of $f(\C^2)$ are
$D_1^2(f,(1,1))=D^2(f)$ (the set of double points),
$f(D_1^1(f,(2)))$ (the set of cross-cap  points) and
$D^3(f,(1,1,1))=D^3(f)$ (the set of  triple points).}
\end{ex}

  The geometric significance of $D^\ell(f,{{\cal P}})$ is given in
  Lemma 2.7  \cite{marar} by W.Marar and D. Mond.
  Given a partition ${{\cal P}} =(r_1,...,r_{\ell})$ of $k$, define
  the projections $\pi_i({{\cal P}}):\C^{n-1+k}\to \C^n,$ for  $1\leq
  i\leq m,$ by $\pi_i({{\cal
  P}})(x,z_1,...,z_k)=(x,z_{r_1+...+r_{i-1}+1}).$

  \begin{lem} (\cite{marar})\label{lema}
  Let ${{\cal P}} =(r_1,...,r_{\ell})$ be a partition of $k$. At a
  generic point (x,z) of $\Delta({{\cal P}})$ we have
  $$
\begin{array}{lll}
{\cal J}_{\Delta}(f,{{\cal P}})=\!\!&{{\cal I}}({{\cal P}})+
(\{\frac{\partial^s f_j}{\partial z^s}\circ \pi_i({{\cal P}})|
j=1,...,p-n+1,1\leq s\leq r_i-1,1\leq i \leq \ell \})\\
&+(\{f_j\circ \pi_1({{\cal P}})-f_j\circ \pi_i({{\cal
P}})|j=1,....p-n+1, 2\leq i\leq \ell\}).
\end{array}$$
  in ${\cal
  O}_{{n-1+k},(x,z)}.$
  \end{lem}

In view of Lemma \ref{lema}, a generic point of $V({\cal
J}_{\Delta}(f,{\cal P}))$ is of the form\\ $({\bf
x},z^1,\cdots,z^1,z^2,\cdots,z^2,\cdots,z^{\ell},\cdots,z^{\ell})$
with ${\bf x}\in \C^{n-1}$, $z^i\in\C$, $z^i$ iterated $r_i+1$
times, $z^i\neq z^j$ for $i\neq j$,
$f(x,z^1)=\cdots=f(x,z^{\ell})$. The local algebra of $f$ at
$({\bf x},z^i)$ is isomorphic to $\frac{\C[z]}{(z^{r_i})}.$

  In  Corollary 2.15 in \cite{marar} is obtained the following result. If
  $f$ is ${\cal A}$-finite then for each partition ${{\cal P}}
  =(r_1,...,r_{\ell})$ of $k$ satisfying $p-k(p-n+1)+\ell \geq 0$, the germ
  of $D^k(f,{{\cal P}})$ at $0$ is either an ICIS of dimension
  $p-k(p-n+1)+\ell$ or is empty.

  \begin{prop}\label{pro3}
Let $F=(u,\overline{f})$ be a versal unfolding of an ${\cal
A}$-finite germ of corank 1. Then for each partition ${\cal P}$ of
$k$ where $p-k(p-n+1)+\ell=0,$ we have the following.
\begin{enumerate}
\item
$D^{\ell}(F,{\cal P})$ is smooth of dimension $s$ or is empty.
\item
${\cal J}_{\Delta}(f,{\cal P})$ is an ICIS.
\item
Let $j_{\ell}:\C^{n-1}\times\C^{\ell}\to \C^{n-1}\times\C^{k}$ be
the embedding with image $\Delta(\cal P)$. Then the surjection
$j_{\ell}^*:{\cal O}_{n-1+k}\to {\cal O}_{n-1+\ell}$ satisfies
$j_{\ell}^*({\cal J}_{\Delta}(f,{\cal P}))={I}(f,{\cal P})$ and
consequently induces an isomorphism
$$j_{\ell}^*:\frac{{\cal O}_{n-1+k}}
{{\cal J}_{\Delta}(f,{\cal P})}\to \frac{{\cal O}_{n-1+\ell}}
{{I}(f,{\cal P})}.$$
\end{enumerate}
\end{prop}
\begin{prf} The items  ${1}$ and ${2}$ follow from \cite{marar}.

{3.} It follows from  Lemma 2.7 \cite{marar} (see also Lemma
\ref{lema} above) that at generic points of $\Delta(\cal P)$ one
has,
$$
\begin{array}{lll}
{\cal J}_{\Delta}(f,{{\cal P}})=\!\!&{{\cal I}}({{\cal P}})+
(\{\frac{\partial^s f_j}{\partial z^s}\circ \pi_i({{\cal P}})|
j=1,...,p-n+1,1\leq s\leq r_i-1,1\leq i \leq \ell \})\\
&+(\{f_j\circ \pi_1({{\cal P}})-f_j\circ \pi_i({{\cal
P}})|j=1,....p-n+1, 2\leq i\leq \ell\}).
\end{array}$$

So generically $j^*_{\ell}({\cal J}_{\Delta}(f,{\cal
P}))={I}(f,{\cal P})$ and  as the two are reduced complete
intersection ideals  coincide we have,
$$\frac{{\cal O}_{n-1+k}}
{{\cal J}_{\Delta}(f,{\cal P})}\cong \frac{{\cal O}_{n-1+\ell}}
{{I}(f,{\cal P})}.$$
\end{prf}

\section{The weighted homogeneous case}

  In what follow we consider weighted homogeneous germs $f:\C^n,0\to\C^p,0.$
  and write $f=(f_1,f_2,...,f_p)$
  The germ $f$ is
  weighted  homogeneous if there exist positive integers $w_1,w_2,....w_n,$ (the
  weights) and  positive integers
  $d_1,d_2,....d_{p}$ (the degrees) such that  for each $f_i$ we have
  $$f_i(t^{w_1}x_1,...,t^{w_n}x_n)=t^{d_i}f_i(x_1,...,x_n),$$
  or equivalently
  $$\Sigma_{j=1}^n w_
  j\alpha_j=d_i,$$ for each monomial
  $x^{\alpha_1}_{1}...x_n^{\alpha_n}$ of $f_i.$

  We give below a formula for calculating the
  number of isolate singularities in the
  case of a weighted homogeneous germ of corank 1.

\begin{theo} \label{for} Lef $f=({\bf
x},f_1,...,f_{p-n+1}):\C^n,0\to\C^p,0$ be an ${\cal A}-$finite
germ de corank 1, with weighted $w_i$ and degree $d_i$ of $f_i$.
Then
$$
\dim_{\C}\frac{{\cal O}_{n-1+\ell}}{{I}(f,{\cal P})}=
\frac{1}{w.w_n^{\ell}}\prod_{j=1}^{p-n+1}
(\prod_{i=1}^{r_i-1}(d_j-iw_n))\prod_{m=2}^{\ell}(\prod_{j=1}^{p-n+1}(d_j-m
w_n))
$$
where $w=w_1.w_2...w_{n-1.}$
\end{theo}
\begin{prf}
According to Proposition \ref{pro3} (3) and Lemma \ref{lema} it
 is enough to compute
 $\dim\frac{{\cal O}_{n-1+k}}{{\cal J}_{\Delta}(f,{\cal P})}.$
 When $f$ is weighted
 homogeneous this dimension can be computed  using  Bezout's theorem
 since ${{\cal J}_{\Delta}(f,{\cal P})}$ is a complete intersection.

The generators of ${{\cal J}_{\Delta}(f,{\cal P})}$ are
$h_{ij}=\frac{\partial^s f_j}{\partial z^s}\circ \pi_i({{\cal
P}})$, $z_i-z_{i+1}$ and $g_{ij}=f_j\circ \pi_1({{\cal
P}})-f_j\circ \pi_i({{\cal P}})$ for each $j=1,...,n-p+1$,
$i=1,...,\ell$ and $s=1,...,r_i-1$. The  degree of  $h_{ij}$ is
$d_j-sw_n$, the degree of $z_i-z_{i+1}$ is $w_n$ for all $i$ and
the degree of $g_{ij}$ is $d_j-mw_n$, where ${m}=2,..\ell$.
 The
product of all the degrees of the generators is therefore
$$\prod^{n-p+1}_{j=1}(\prod^{r_i-1}_{i=1}(d_j-iw_n))\prod_{m=2}^{\ell}
(\prod_{j=1}^{p-n+1}(d_j-m w_n)).w_n^{k-\ell}.$$ Since ${{\cal
J}_{\Delta}(f,{\cal P})}$ is a weighted homogeneous complete
intersection, we can apply Bezout's Theorem, hence its colength is
given by

$$
\frac{1}{w_n^{k}.w}\prod_{j=1}^{p-n+1}
(\prod_{i=1}^{r_i-1}(d_j-iw_n))\prod_{m=2}^{\ell}(\prod_{j=1}^{p-n+1}(d_j-m
w_n)).w_n^{k-\ell}
$$
as required.

\end{prf}

\begin{ex} {\rm Let $f:\C^2,0\to\C^3,0$ be an ${\cal A}$-finite germ
given by $f(x,y)=(x,xy+y^3,y^4)$ (see \cite{mond}). Then the
partition ${\cal P}$ is $(1,1,1)$ and $(2)$ and the weights and
degrees are $d_1=3$, $d_2=4$ and $w_1=2$, $w_2=1.$ Therefore
$\sharp{\cal Q }(f,{(1,1,1)})=1$, $\sharp{\cal Q }(f,{(2)})=2$,
that is $f$ has 1 triple point and  2 cross caps (see
\cite{mond}).

Let $f:\C^3,0\to\C^4,0$ be an ${\cal A}$-finite germ  given by
$f(x,y,z)=(x,y,yz+z^4,xz+z^3).$ The partition ${\cal P}$ is
$(1,1,1)$ and $(2,1)$ and the weights and degrees are $d_1=4$,
$d_2=3$ and $w_1=2$, $w_2=3$, $w_3=1.$ Therefore $$\sharp{\cal Q
}(f,{(2,1)})=\frac{(4-3)(3-1)(4-1)(3-1)}{2.3.1^2}=2$$  and
$\sharp{\cal Q }(f,{(1,1,1,1)})=0$, that is $f$  has 2
singularities of type ${\cal Q }(f,{(2,1)})$ and  has not
quadruple points. In particular applying Theorem  \ref{for} or the
method of Example \ref{exe3},  all  corank-1 simple germs
$f:\C^3,0\to\C^4,0$ classified by Houston and Kirk in \cite{hk}
satisfy $$\sharp{\cal Q }(f,{(1,1,1,1)})=0.$$

Observe that if ${\cal P}=(2,1)$ and $f:\C^3,0\to\C^4,0$, $V({\cal
J}_{\Delta}(f,(2,1)))\subset \C^5$ and $V(I(f,(2,1)))\subset \C^4$
are isomorphic. Hence $\mu(V({\cal
J}_{\Delta}(f,(2,1))))=\mu(V(I^3(f)+{\cal I}(2,1)))=\mu(D^3(f)|H)$
where $H=V({\cal I}(2,1))$, therefore as $V({\cal
J}_{\Delta}(f,(2,1)))$ is ICIS zero-dimensional, we have
$$\mu(D^3(f)|H)=\dim_{\C}\frac{{\cal O}_5}{(I^3(f)+{\cal I}(2,1))}-1=
\dim_{\C}\frac{{\cal O}_4}{I(f,(2,1))}-1=\sharp{\cal Q
}(f,{(2,1)})-1,$$ as Houston and Kirk compute $\mu(D^3(f)|H)$, we
calculate $\sharp{\cal Q }(f,{(2,1)})$ in the table below.

$
\begin{array}{llcc}
 Label&Singularity&\mu(D^3(f)|H)&\sharp{\cal Q }(f,{(2,1)})\\
P_1&(x,y,yz+z^4,xz+z^3)&1&2\\
 P_2&(x,y,yz+z^5,xz+z^3)&2&3\\
  P_3^k&(x,y,yz+z^6+z^{3k+2},xz+z^3)&3&4\\
  P_4^1&(x,y,yz+z^7+z^8,xz+z^3)&4&5\\
P_4&(x,y,yz+z^7,xz+z^3)&4&5\\
Q_k&(x,y,xz+yz^2,y^kz+z^3)&1&2\\
R_k&(x,y,xz+z^3,yz^2+z^4+z^{2k-1})&2&3\\
S_{j,k}&(x,y,xz+y^2z^2+z^{3j+2},z^3+y^kz)&3&4
\end{array}
$}
\end{ex}

\section{ Necessary and sufficient conditions
for${\cal A}$-finiteness of map-germs} In this section we define
new numerical invariants associated to each partition ${\cal P}$
of $k$ with $p-k(p-n+1)+\ell\geq 0$. We show  that a germ $f$ is
${\cal A}$-finite if and only if these invariants are finite.

For any partition ${\cal P}$ of $k$, $p-k(p-n+1)+\ell\geq 0,$
denote by $H_{\cal P}$ the map defined by
$$H_{\cal P}:\C^{n-1+\ell}\to\C^{(p-n+1)(k-\ell)}$$ with
components  $\frac{\partial^s
  {f}_j}{\partial z^s}\circ \pi_i({{\cal P}})$ for $j=1,...,p-n+1,1\leq
  s\leq r_i-1,1\leq i \leq \ell$  and

$$G_{\cal P}:\C^{n-1+\ell}\to\C^{(p-n+1)(\ell-1)}$$
 with components
${f}_j\circ \pi_1({{\cal P}})-{f}_j\circ
  \pi_i({{\cal P}})$ for $j=1,....p-n+1, 2\leq i\leq \ell.$

 Let
 $$F_{\cal P}=(G_{\cal P},H_{\cal P}):\C^{n-1+\ell}\to
 \C^{(p-n+1)(k-1)},$$

and define the following number
 $$N(f,{\cal P})=
 \dim_{\C}\frac{{\cal O}_{n-1+\ell}}{F^*_{\cal P}
 ({\cal M}_{(p-n+1)(k-1)})+J(F_{\cal P})}$$

 where $J(F_{\cal P})$ denotes the ideal
 of  $(p-n+1)(k-1)\times(p-n+1)(k-1)$ minors
 and ${\cal M}_{(p-n+1)(k-1)}$ is the maximal ideal in ${\cal
 O}_{(p-n+1)(k-1)}$.

 We first show that $N(f,{\cal P})$ is
${\cal A}$-invariant.
 \begin{prop}
 Let $f_1$ and $f_2$ be corank 1, ${\cal A}$-finite
   and  ${\cal A}$-equivalent maps germs. Then $N(f_1,{\cal P})=N(f_2,{\cal P}).$
 \end{prop}
 \begin{prf} It is  sufficient to proof that if $f_1$ is ${\cal
 A}$-equivalent to $f_2$, then the germs defined above, $F_{1{\cal
 P}}$ associated to $f_1$ and $F_{2{\cal P}}$ associated to $f_2$, are ${\cal
 K}$-equivalent.
 The last statement follows from \cite{mond}.
 \end{prf}

\begin{theo} Let  $f:\C^n,0\to\C^p,0$, $n<p$ be a corank 1 germ.
Then following statements are equivalent:
\begin{enumerate}
\item $f$ is an ${\cal A}$-finite.
\item $N(f,{\cal P})<\infty$, where ${\cal P}$ is
$(1,1),...,(1,...,1)$  $k$-times and $\sharp{\cal Q}(f,{\cal
P})<\infty$ for each  partition  ${\cal P}$  of $k$ where
$p-k(p-n+1)+\ell\geq 0$.
\end{enumerate}
\end{theo}
\begin{prf}
$(1)\Longrightarrow (2)$  $f$ is ${\cal A}$-finite if and only if
for any representative of $f$ there exist  neighborhoods $U$ of
$0$ in $\C^n$ and $V$ of $0$ in $\C^p$, with $f(U)\subset V$, such
that for all $y\neq 0$ in $V$, the multi germ
 $f:(U,f^{-1}(y)\cap \Sigma(f))\to (V,y)$ is stable. If follows from
 Proposition 2.13 in \cite{marar} that if $f$ is
 of type $\Sigma^{1_{r_i},0}$ in $({\bf x},y_i),$ with
 $y_i\in f^{-1}(y)\cap\Sigma(f)$ then $F_{\cal P}$
 defining  $V(J_{\Delta}(f,{\cal P}))$ is a
 submersion at $({\bf x},y_1,...,y_k)$. Then at every point of
 $V(J_{\Delta}(f,{\cal P}))$, distinct from 0, the $(p-n+1)(k-1)$ functions
 generating $J_{\Delta}(f,{\cal P})$ define a submersion,
 and so $V(J_{\Delta}(f,{\cal P}))$ is ICIS. Therefore $V(F^*_{\cal P}({\cal
M}_{(p-n+1)(k-1)})+J(F_{\cal P})),$ for each partition $k$ where
$p-k(p-n+1)+\ell\geq 0$, is zero in $\C^{n-1+\ell}$. Hence
$V(F^*_{\cal P}({\cal M}_{(p-n+1)(k-1)})+J(F_{\cal P}))\subset
\{0\}$ and by Nullstellensatz we have $N(f,{\cal P})<\infty.$

$(1)\Longleftarrow (2)$ if $N(f,{\cal P})<\infty$ we have
$V(F^*_{\cal P}({\cal M}_{(p-n+1)(k-1)})+J(F_{\cal
P}))\subset\{0\}$, that is  $V(F^*_{\cal P}({\cal
M}_{(p-n+1)(k-1)}))$ have isolated singularities in $\{0\}$ for
each partition ${\cal P}$, in particular for
$(1,1),...,(1,...,1)$-$k-1$ times and for all partition ${\cal P}$
of $k$ where, $p-k(p-n+1)+\ell=0.$ We choose  representative
$f:U\subset \C^n\to V\subset \C^p$, such that the representatives
induced  by the germs of $F^{-1}_{\cal P}(0)$ are differentiable.
With this we have the multi germ $f:(U,S)\to(V,z),$ where
$S=f^{-1}(z),$ which is stable for $z\neq 0$.  The result follows
then by the geometric criterion of Mather and Gaffney \cite{cg}.
\end{prf}
\begin{ex} {\rm Let $f:(\C^3,0)\to(\C^4,0)$ be given by
$$f(x,y,z)=(x,y,xz+z^3,yz^2+g(z))$$ where
$g(z)=\Sigma_{i=1}^ka_iz^i.$ Then $f$ is ${\cal A}$-finite  for
all $a_i$, The partitions ${\cal P}$ are $(1,1,1,1)$, $(1,1,1)$,
$(1,1)$, $(2,1)$ and $(2)$ and the stable types of $f$ are:
$D^4(f,(1,1,1,1))$ which is empty, the triple points curve
$D^3(f,(1,1,1))$, the types 2-dimensional $D^2(f,(1,1))$, the
types 1-dimensional $D^1(f,(2))$
 and the zero-dimensional types $D^2(f,(2,1))$. We compute the  maps $F_{\cal P }$ for all
partitions ${\cal P}$ using Maple.

If ${\cal P}=(1,1,1)$, then $F_{\cal
P}(x,y,z,v,w)=(x-(zv+zw+vw),z+v+w,V_1^3(g),y+V_2^3(g))$ where
$V_i^k(g)$ is  computed  using the Definition given in the section
\ref{sec4}.

If ${\cal P}=(1,1)$, then $F_{\cal
P}(x,y,z,v)=(x+z^2+zv+v^2,y(z+v)+V_1^2(g))$.

 Using the software  Singular \cite{singular} we have
$$
\begin{array}{lll}
N(f,(1,1,1))&=&\dim_{\C}\frac{{\cal O}_5
}{(x-(zv+zw+vw),z+v+w,V_1^3(g),y+V_2^3(g),J(F_{\cal P}))}\\
 &=&\left\{
\begin{array}{lll}
0 & \mbox{if}\,\, a_1=a_2=0\\
2& \mbox{if}\,\,{\ell}(g)=3\\
5& \mbox{if}\,\,{\ell}(g)=4\\
20& \mbox{if}\,\,{\ell}(g)=5\\
17&\mbox{if}\,\,{\ell}(g)=6\\
26&\mbox{if}\,\,{\ell}(g)=7
\end{array}
\right. %BIBLIOTECARIA777@BOL.COM.BR, paulinha, (19) 9135-9922
\end{array}
$$
where ${\ell}(g)$ denote the lowest degree in $g$.
$$ N(f,(1,1))=\left\{
\begin{array}{l}
{\ell}(g)-1\,\, \mbox{if}\,\,{\ell}(g)\,\,\mbox{is odd}\\
{\ell}(g)\,\,\,\mbox{if}\,\, {\ell}(g) \,\,\mbox{is even.}
\end{array}
\right.
$$

The ideal that defines $D^2(f,(2,1))$ is given  by ${\cal
I}^3(f)+{\cal I}({\cal P})$ where ${\cal I}^3(f)$ defines
$D^3(f)=D^3(f,(1,1,1))$ and ${\cal I}({\cal P})=v-w$. Then
$D^2(f,(2,1))=D^3(f,(1,1,1))\cap  H$, and
$D^2(f,(2,1))=V(x-(zv+zw+vw),z+v+w,V_1^3(g),y+V_2^3(g),v-w)\subset
\C^{2}\times\C^3$. If ${\cal P}=(2,1)$, we have $F_{\cal
P}(x,y,z,v)=(x-(zv+zv+v^2),z+2v,V_1^3(g)(z,v),y+V_2^3(g)(z,v))$,
then $N(f,(2,1))<\infty$.

The ideal ${\cal I}^2(f)+{\cal I}({\cal P})$  defines $D^1(f,(2))$
where ${\cal I}^2(f)$ defines $D^2(f)=D^2(f,(1,1))$ and ${\cal
I}({\cal P})=v-z$. We have $D^1(f,(2))=D^2(f,(1,1))\cap H$, where
$H=V(v-w)\subset \C^2\times\C^2$ is a hyperplane. Therefore
$D^1(f,(2))= V(x+z^2+zv+v^2,y(z+v)+V_1^2(g),v-z)\subset
\C^{2}\times\C^2$. If ${\cal P}=(2)$ we have $F_{\cal
P}(x,y,z)=(x+3z^2,2zy+V_1^2(g)(z,z))$, and therefore
$N(f,(2))<\infty$ for all $a_i$.}
\end{ex}

\begin{rem}{\rm Observe that if $g(z)=z^2+z^7$ and $g(z)=z^5+z^6+z^7$, then
the germs
$$f_1(x,y,z)=(x,y,xz+z^3,yz^2+z^2+z^7)$$ $$f_2(x,y,z)=
(x,y,xz+z^3,yz^2+z^5+z^6+z^7)$$  are not ${\cal A}$-equivalent as
$N(f_1,(1,1))\neq N(f_2,(1,1))$.}
\end{rem}

The number of invariants
 above involved depends on the dimensions $(n,p)$, and this number is large when $n$ and $p$
 are large. It is then natural to ask:
Fixing a pair $(n,p)$, what is the minimum number of invariants
$N(f,{\cal P})$ are necessary and sufficient to ensure ${\cal
A}$-finiteness of the germ $f$?, Then we have the following.

{\bf Problem:} {\it Fixing a pair $(n,p)$, The numbers of
invariants $N(f,{\cal P})$ where ${\cal P}$ is the partition of
$k$ such that $p-k(p-n+1)+\ell=0$,  are necessary and sufficient
to ensure ${\cal A}$-finiteness of the germ $f$?}

This question has been answered for the cases:
\begin{rem}{\rm when $n=2$ and $p=3$,
the answer is given by  D. Mond in \cite{mond}. Then $f$ is ${\cal
A}$-finite if and only if the number of cross caps $C(f)$, the
number triple points $T(f)$  and $N(f)<\infty$ are finites. It
turns our  that $N(f)=N(f,(1,1))$.

 In the case $n=p=2$ Gaffney and Mond \cite{GM1} showed  that $f$ is
${\cal A}$-finite if and only  if the number of cusps number
$c(f)$, the number double folds points $d(f)$ are finite.}
\end{rem}

% este inavariante es facil calcular usando o singular.

\medskip
\noindent{\bf Acknowledgements.} I would like to thanks to Farid
Tari for discussions and numerous remark on mathematics and
English usage of the paper.

\medskip

\noindent {\it Address:} Instituto de Ci\^encias Matem\'atica e de
Computa\c c\~ao ICMC-USP\\ Av. do Trabalhador  S\~ao-Carlense, 400
- Centro - Cx. Postal 668 \\ S\~ao Carlos - S\~ao Paulo - Brasil
CEP 13560-970
\\ \noindent E-mail:
vhjperez@icmc.usp.br

\end{document}